\documentclass[a4paper,12pt]{article}

\usepackage{pstricks}
\usepackage{graphicx,psfrag}
\usepackage{amsmath}
\usepackage{amssymb}
\usepackage{amsthm}
\usepackage{color}

      \newcommand {\al}   {\alpha}          \newcommand {\bt}  {\beta}
                
      \newcommand {\del}  {\delta}          
              \newcommand {\ve}   {\varepsilon}
      \newcommand {\z }   {\zeta}           \newcommand {\vphi} {\varphi}
      \newcommand {\lam}  {\lambda}

                \newcommand {\Om}  {\Omega}
      \newcommand {\pl}   {\partial}

           \newcommand {\UUU}  {{\cal U}}
             
      \newcommand {\RRR}  {{\mathbb R}}

     \newcommand {\beq}  {\begin{equation}}  \newcommand {\beqo}  {\begin{equation*}}
      \newcommand {\eeq}  {\end{equation}}   \newcommand {\eeqo}  {\end{equation*}}

      \newtheorem{theorem}{Theorem}

      \newtheorem{cor}{Corollary}

\newcommand{\re}{{\mathbb R}}

\newtheorem{remark}{Remark}

\title{Local minima in Newton's aerodynamical problem  \\
and inequalities between norms of partial derivatives}

\author{Alexander Plakhov\thanks{Center for R{\&}D in Mathematics and Applications, Department of Mathematics, University of Aveiro, Portugal and Institute for Information Transmission Problems, Moscow, Russia, plakhov@ua.pt} \and Vladimir Protasov\thanks{DISIM, University of L'Aquila, Italy and Moscow State University, Russia}}

\begin{document}
\maketitle

\begin{abstract}
The problem considered first by I.\,Newton (1687)  consists in finding a surface of the minimal frontal resistance
 in a  parallel flow of non-interacting point particles. The standard formulation assumes that
 the surface is convex with a given convex base~$\Omega$ and a bounded altitude.
 Newton found the solution for surfaces of revolution.  Without this assumption
 the problem is still unsolved, although many important results have been obtained  in the last decades.
 We consider the problem to characterize the domains~$\Omega$ for which the flat surface
 gives a local minimum.   We show that this problem can be reduced to
 an inequality between $L_2$-norms of partial derivatives for bivariate concave
 functions on a convex domain that vanish on the boundary. Can the ratio between
  those norms be arbitrarily large? The answer depends on the geometry of the domain.
 A complete criterion is derived, which also solves the local minimality problem.
\end{abstract}

\begin{quote}
\end{quote}

\begin{quote}
\end{quote}

\section{Introduction}

\subsection{Newton's problem of minimal resistance}

We answer one question in Newton's  minimal resistance problem by means of
an inequality between norms of partial derivatives. The inequality deals with
$L_2$-norms of derivatives~$u_x$ and~$u_y$ of a concave bivariate function~$u(x,y)$
which is given on a convex domain and vanishes on its boundary.
We derive a geometrical criterion on the domain that ensures that
the ratio~$\|u_x\|_{2}/\|u_y\|_2$ can be arbitrarily large. Then using this result we
decide when a flat surface gives a local minimum in the Newton problem.
We begin with notation and statements of the problems.

\begin{quote}
Given a continuous function $f : \RRR^2 \to \RRR$, a bounded convex set with nonempty interior $\Om$, and $M > 0$, 
\beq\label{FunRes}
 \text{minimize} \qquad F(u) = \int_\Om f(\nabla u(x,y))\, dx\, dy
\eeq
in the class of continuous concave functions $u : \Om \to \RRR$ satisfying $0 \le u \le M$.
\end{quote}
This problem goes back to Isaac Newton and has a simple mechanical interpretation.

Namely, consider a convex body in Euclidean space $\RRR^3$ with the coordinates $x$, $y$, $z$ and a homogeneous flow 
composed of point particles moving downward with the velocity $(0,0,-1)$. The flow is so rare that the particles do not interact with each other; they hit the body's surface and are reflected in accordance with a certain {\it scattering law}, and then go away. The scattering law can be understood as a family of probability measures in $S^2 \times S^2$ (the set of pairs (velocity of incidence, velocity after reflection)\,) parameterized by vectors from $S^2$ (outward normal to the body at the point of reflection). As a result of particle-body interaction, the drag force is created.

Let the front (exposed to the flow) part of the body's boundary be the graph of a concave function $u$, with the domain $\Om$ of the function being the projection of the body on the $xy$-plane; then the $z$-component of the drag force (usually called the {\it resistance}) equals $2\rho F(u)$, where $F$ is defined according to \eqref{FunRes},\, $\rho$ is the density of the flow, and the function $f$ is determined by the scattering law. Thus, problem \eqref{FunRes} is to find the body having the minimal resistance among all convex bodies with the fixed projection $\Om$ on the $xy$-plane and with the fixed projection $[0,\, M]$ on the $z$-axis.

In particular, if the particles are reflected in the perfectly elastic (billiard) manner, $f$ is equal to
\beq\label{fu}
f(\xi_1,\xi_2) = \frac{1}{1 + \xi_1^2 + \xi_2^2},
\eeq
and the (normalized) resistance equals
\beq\label{FunBil}
 F(u) = \int_\Om \frac{1}{1+|\nabla u(x,y)|^2}\, dx\, dy.
\eeq

\begin{figure}[h]
\centering
\includegraphics[scale=0.3]{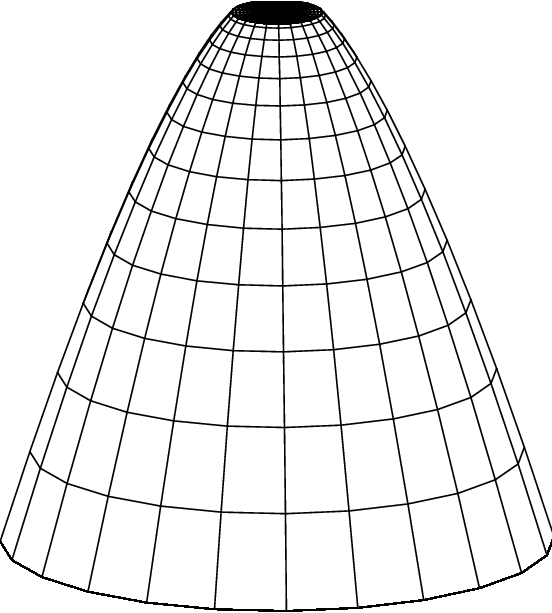}
\caption{A solution of Newton's problem.}
\label{figNewtonOpt}
\end{figure}

Newton himself considered the problem of least resistance, in the case of elastic reflections, in the narrower class of convex {\it rotationally symmetric} bodies. Then the set $\Om$ is a circle with the radius, say, equal to $L$, the function $u$ is radially symmetric, $u(x,y) = \vphi(\sqrt{x^2+y^2})$, and the functional \eqref{FunBil} takes the form
$2\pi \int_0^L \frac{1}{1 + \vphi'^2(r)}\, r\, dr$. Thus, the problem amounts to minimizing the integral
$$
\int_0^L \frac{1}{1 + \vphi'^2(r)}\, r\, dr
$$
in the class of concave monotone non-increasing functions $\vphi: [0,\, L] \to \RRR$ satisfying $0 \le \vphi \le M$. Newton provided the solution in geometric terms and without proof. Here is how it looks when $M/L \approx 2$ (see Fig.~\ref{figNewtonOpt}).

The interest to Newton's problem revived in the early 1990s starting from the seminal paper \cite{BK} by Buttazzo and Kawohl. They posed the problem of minimizing \eqref{FunBil}, without the symmetry condition, and conjectured that the solution does not coincide with Newton's one. Since then a significant work on the problem has been made (see, e.g., \cite{BrFK,BFK,MMO,P,boundary,develop,LZ}). Some of the obtained results hold true also for the more general problem \eqref{FunRes}. It was proved that the solution $u$ always exists \cite{BK} and does not coincide with Newton's one \cite{BrFK}; that is, by relaxing the symmetry condition one can find bodies that are more streamlined than Newton's one. Further, $|\nabla u| \not\in (0,\, 1)$ \cite{BFK}, and if $|\nabla u|$ takes the zero value then it takes values arbitrarily close to 1 \cite{MMO}. Additionally, it was proved that
any optimal body (that is, the set $(x,y) \in \Om$,\, $0 \le z \le u(x,y)$ with $u$ optimal) is the convex hull of the closure of the set of its singular points \cite{P}.

In \cite{develop} and \cite{LZ} the problem is studied in several narrower classes of convex bodies. The bodies in these classes are convex hulls of two curves, where the former curve is the fixed circumference $\pl\Om \times \{0\}$ and the latter curve is to be optimized. In \cite{develop} the latter curve is a convex closed curve in the horizontal plane $z=M$, while in \cite{LZ} it is a convex curve in a vertical plane situated between the planes $z=0$ and $z=M$ with the endpoints at two opposite points of the circumference $\pl\Om \times \{0\}$.

Sometimes the zero function $u \equiv 0$ is a solution of problem \eqref{FunRes}. For example, the maximum of Newtonian resistance when $f$ is given by \eqref{fu} (or, equivalently, the minimum of $F$ with the integrand equal to $-f$) is attained at the zero function. Besides, it is proved in \cite{boundary} that, under quite mild conditions on $\Om$ and $f$, the solution $u$ of problem \eqref{FunRes} satisfies $u\rfloor_{\pl\Om} = 0$. In the present paper we are interested in the question, if the zero function is a local minimum of the problem. (Note that in this setting the value of $M$ does not make difference.) In the next section \ref{ss12} we formulate the minimization problem, in a slightly modified form as compared with the problem of minimal resistance, and state the main results, Theorems \ref{t1} and \ref{t2}. These theorems are proved in sections \ref{st1a}, \ref{st1b}, and \ref{st2}.

\subsection{Formulation of the problem and main results}\label{ss12}

Consider a convex bounded set with nonempty interior $\Om$ on the plane $\RRR^2$ with the coordinates $x,\, y$. Denote by $\UUU_\Om$ the set of continuous concave functions $u : \Om \to \RRR$ with bounded gradient such that $u\rfloor_{\pl\Om} = 0$.

Consider a tangent cone to $\Om$ with the vertex at a singular point of $\pl\Om$. A straight line that intersects the tangent cone through its vertex is called an {\it angular} support line.

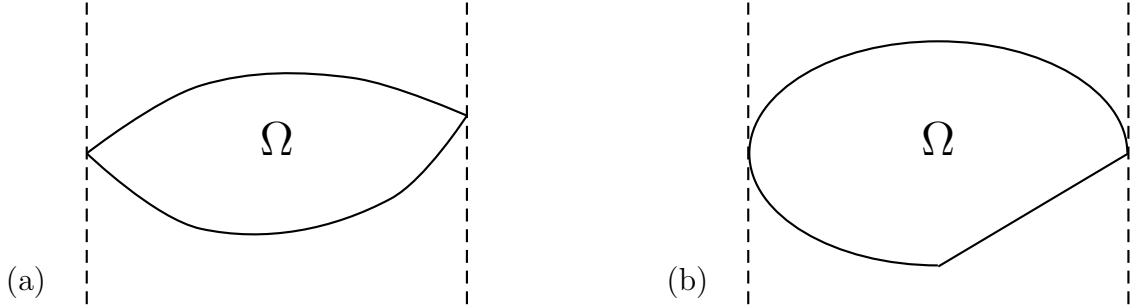
\begin{figure}[h]
\begin{picture}(0,130)
 \scalebox{1}{

 \rput(1.3,2.1){
 \psline[linestyle=dashed](0,-2)(0,2) \psline[linestyle=dashed](5,-2)(5,2)
 \pscurve(0,0)(1.5,0.9)(3.5,1)(5,0.5) \pscurve(0,0)(1.5,-1)(4,-0.6)(5,0.5)
\rput(2.5,0.2){\scalebox{1.5}{$\Om$}} \rput(-0.8,-1.7){(a)}
}

\rput(10,2.1){
\psellipse(2.5,0)(2.5,1.5)
\pspolygon[fillstyle=solid,fillcolor=white,linewidth=0pt,linecolor=white](2.5,0)(5,0)(5,-1.5)(2.5,-1.5)
\psline(5,0)(2.5,-1.5)
\psline[linestyle=dashed](0,-2)(0,2) \psline[linestyle=dashed](5,-2)(5,2)
\rput(2.5,0.2){\scalebox{1.5}{$\Om$}} \rput(-0.8,-1.7){(b)}
}

}
\end{picture}
\caption{In figure (a), both vertical lines of support to $\Om$ are angular, while in figure (b), both lines are not angular: the left line is tangent, and the right one is half-tangent to $\Om$.}
\label{fig:ang}
\end{figure}

In Fig~\ref{fig:ang}, examples illustrating the notion of angular line are given.

The main results of this article are formulated in the following Theorems \ref{t1} and \ref{t2}.

\begin{theorem}\label{t1}
(a) If at least one vertical (that is, parallel to the $y$-axis) line of support is not angular then
$$
\sup_{\stackrel{u\in\UUU_\Om}{u\not\equiv0}}\, \frac{\int_\Om u_x^2\, dx\, dy}{\int_\Om u_y^2\, dx\, dy} = \infty.
$$

(b) If both vertical lines of support are angular then
$$
\sup_{\stackrel{u\in\UUU_\Om}{u\not\equiv0}}\, \frac{\int_\Om u_x^2\, dx\, dy}{\int_\Om u_y^2\, dx\, dy} =: K < \infty.
$$
\end{theorem}
\begin{remark}\label{r10}
{\em Note that without the concavity condition on the function~$u$ in Theorem~\ref{t1},
the supremum of the ratio~$\dfrac{\int_\Om u_x^2\, dx\, dy}{\int_\Om u_y^2\, dx\, dy}$
is infinite for every convex domain~$\Omega$. To see this it suffices to
take an arbitrary nonzero~$C^1$ function~$u_0(x,y)$ that vanishes on~$\, \partial \Omega$
and multiply it by a highly oscillating function~$\varphi(x)$, for example,
by~$\varphi_N(x) = \sin \,Nx$. Then
the ratio of norms of partial derivatives for the function~$u_0(x,y)\varphi_N(x)$
tends to infinity as~$N\to \infty$. }
\end{remark}

\begin{remark}\label{r20}
{\em  Theorem~\ref{t1} can be attributed to the so-called inequalities between derivatives,
which is an important branch of the approximation theory known from the classical
results of Bernstein, Kolmogorov, Landau, Nikolsky, Calderon, Zygmund, etc. Such inequalities involve
norms of derivatives of various orders in different functional spaces.
 See~\cite{HR} -- \cite{M} for an extensive bibliography. We, however, have not meet
special  inequalities for concave functions on convex domains.  To the best of our knowledge, the result  of~Theorem~\ref{t1} is new.
}
\end{remark}

Now we apply Theorem~\ref{t1} to the following problem of minimal resistance: given a continuous function $f : \RRR^2 \to \RRR$, find

\beq\label{problem}
\inf_{u\in\UUU_\Om} F(u), \quad \text{where} \quad F(u) = \int_\Om f(\nabla u(x,y))\, dx\, dy.
\eeq

Note that in general a function $u \in \UUU_\Om$ is not a $C^1$ function; however, it is differentiable almost everywhere. We define the $C^1$ norm in the linear space of continuous and almost everywhere differentiable functions $u$ with the bounded gradient on $\Om$ as follows:
\beq\label{C1norm}
\|u\|_{C^1} = \sup_{(x,y)\in\Om} |u(x,y)| + \sup_{(x,y)\in\Om} |u_x(x,y)| + \sup_{(x,y)\in\Om} |u_y(x,y)|,
\eeq
the suprema in the second and third terms being taken over the points $(x,y)$ where the derivative $u_x$ and $u_y$, respectively, exists.

Here we consider the question, if the zero function $u \equiv 0$ is a local minimum of problem \eqref{problem} in the $C^1$ norm. Theorem \ref{t2} gives the answer to this question.

\begin{theorem}\label{t2}
Let  $f$ be twice differentiable at $(0,0)$.

(a) If the quadratic form $f''(0,0)$ is positive definite then $u \equiv 0$ is a local minimum of problem \eqref{problem}.

(b) If $f''(0,0)$ is negative definite then $u \equiv 0$ is not a local minimum.

(c) Let $f''(0,0)$ be indefinite. Denote by $-a < 0 < b$ its eigenvalues. Here one should consider two cases.

\hspace*{10mm}(i) If at least one line of support parallel to the positive eigendirection (that is, the eigendirection corresponding to the positive eigenvalue) of $f''(0,0)$ is not angular, then $u \equiv 0$ is not a local minimum.

\hspace*{10mm}(ii) If both lines of support parallel to the positive eigendirection are angular then for $b/a < K$, $u \equiv 0$ is not a local minimum and for $b/a > K$, $u \equiv 0$ is a local minimum, where the value $K > 0$ is defined in claim (b) of Theorem \ref{t1}.
\end{theorem}

\begin{remark}\label{remK}
Note that $K$ depends only on the domain $\Om$ and on the positive and negative eigendirections of $f''(0,0)$.
\end{remark}

Claim (c) is of course the most important in Theorem \ref{t2}. The case when $f''(0,0)$ is semidefinite remains open.
\smallskip

The following statement is a direct consequence of~Theorem \ref{t2}

\begin{cor}\label{cor1}
If $f''(0,0)$ is indefinite and $\pl\Om$ contains at most one singular point then $u \equiv 0$ is not a local minimum of problem \eqref{problem}.
\end{cor}

\section{Proof of Theorem \ref{t1}: claim (a)}\label{st1a}

\begin{figure}[h]
\begin{picture}(0,160)
\scalebox{1}{

\rput(7.5,0.1){
\psarc(0,0){5}{0}{180} \psline(-5,0)(5,0) \psline[linestyle=dashed](3.54,3.54)(6.04,1.04)(5,0)
\psdots(4.83,1.29) \psline[linestyle=dotted,linewidth=1.2pt](5.18,0)(4.1,4)  \rput(4.4,4){\scalebox{1}{$l_{\theta}$}}
\psdots[dotsize=2.5pt]
(4.719,0.87)(4.514,1.63) 
\psline(4.719,0.87)(4.514,1.63)
\rput(-0.5,2){\scalebox{3}{$\Om$}}

\rput(5.1,1.3){\scalebox{1}{$\xi$}} \rput(4.52,0.62){\scalebox{1}{$A$}} \rput(4.25,1.7){\scalebox{1}{$B$}}
}}
\end{picture}
\caption{Here $\Om$ is the upper half-circle, and $\Om'$ contains $\Om$ and is bounded by the dashed broken line on the right hand side.}
\label{fig:AB}
\end{figure}
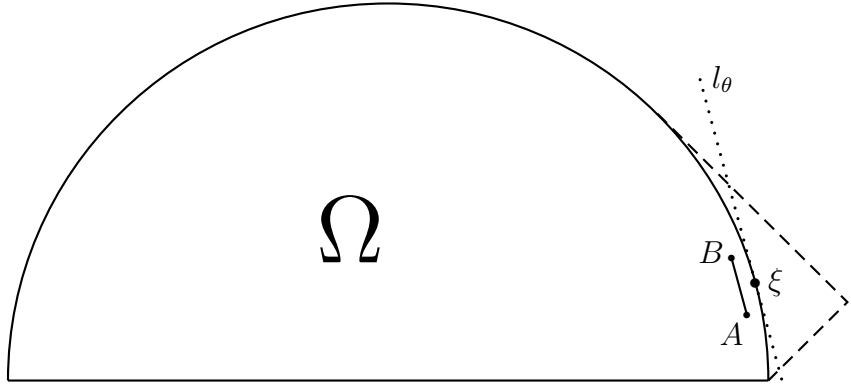

Denote by $l_\al$ the line of support with the outward normal $(\cos\al, \sin\al)$,\, $-\pi < \al \le \pi$, and by $\Pi_\al$ the closed half-plane bounded by $l_\al$ and containing $\Om$. By the hypothesis, one of the lines $l_0$, $l_\pi$ is not angular; let it be $l_0$. Fix a value $0 < \vphi < \pi/2$ and let
$$
\Om' = \Om'_\vphi = \cap_{|\al|\ge\vphi} \Pi_\al;
$$
see Fig.~\ref{fig:AB}. Since $l_0$ is not angular, there is an open arc of $\pl\Om$ contained in the interior of $\Om'$. All lines of support at points of this arc correspond to angles in $(-\vphi,\, \vphi)$.

Choose a regular point $\xi$ on the arc, and let $l_\theta$,\, $-\vphi < \theta < \vphi$, be the (unique) line of support at this point.

Denote $r = \text{dist}(\xi, \pl\Om') > 0$. Take $\ve > 0$ and draw the line $l_{\theta,\ve}$ contained in the half-plane $\Pi_\theta$ (and therefore, parallel to $l_\theta$) at the distance $\ve$ from $l_\theta$. For $\ve$ sufficiently small, the intersection of $l_{\theta,\ve}$ with $\Om$ is a nonempty line segment with the length $\lam(\ve)$ satisfying $\lam(\ve)/\ve \to \infty$ as $\ve \to 0$. Choose two interior points of the segment $A = A_\ve$ and $B = B_\ve$ so as the segment $AB$ belongs to the $(r/2)$-neighborhood of $\xi$ and its length satisfies
\beq\label{AB}
\frac{|A_\ve B_\ve|}{\ve} \to \infty \qquad \text{as} \quad \ve \to 0.
\eeq
See Fig.~\ref{fig:AB}.

Take the convex hull of the sets $\Om \times \{0\}$ and  $AB \times \{1\}$ in $\RRR^3$. The upper part of its boundary is the graph of a concave function $u = u_{\vphi,\ve}$ (the lower part is $\Om \times \{0\}$).

The graph of $u$ is the union of segments with the former endpoint on the boundary of $\Om$ and the latter one on the segment $AB$. Consider the natural projection (under the map $(x,y,z) \mapsto (x,y)$) of any such segment on the $xy$-plane. Let $\Om_1 = \Om_1(\vphi,\ve)$ be the union of projections of segments with the former endpoint on $\pl\Om \cap \pl\Om'$, and $\Om_2 = \Om_2(\vphi,\ve)$ be the union of projections of segments with the former endpoint on $\pl\Om \setminus \pl\Om'$.

Note that $\Om$ is the union of projections of segments with the latter endpoint at $A$ or at $B$ and of two trapezoids (which may degenerate to triangles). One base of each trapezoid coincides with $AB$, and the other one is the intersection of $\pl\Om$ with $l_\theta$ or of $\pl\Om$ with $l_\theta+\pi$ (each one of the intersections can degenerate to a point, the point $\xi$ in the former case.) Note also that the former trapezoid lies in $\Om_2.$

We have $\Om_1 \cup \Om_2 = \Om$ and $\Om_1 \cap \Om_2 = [A, B]$.
Correspondingly,
$$
\int_\Om |\nabla u|^2\, dx dy = I_1 + I_2,
$$
where
$$
I_j  = I_j(\vphi,\ve) = \int_{\Om_j} |\nabla u|^2\, dx dy, \quad j = 1,\, 2.
$$
For regular points of $\Om_1$ one has $|\nabla u| \le \dfrac{1}{r/2}$, hence
$$
I_1 \le |\Om| (2/r)^2.
$$
Here and in what follows, $|\Om|$ means the area of $\Om$.

Further, taking into account that the triangle $AB\xi$ belongs to $\Om_2$, its area equals $\frac 12\, \ve |A B|$, and $|\nabla u| = 1/\ve$ at all points of the triangle, and using \eqref{AB} one obtains
$$
I_2 \ge \int_{\bigtriangleup A_\ve B_\ve \xi} |\nabla u|^2\, dx dy = \frac 12\, \ve |A_\ve B_\ve| \cdot (1/\ve)^2 \to \infty \quad \text{as} \quad \ve \to 0.
$$

On the other hand, at each regular point of $\Om_2$ one has $|u_y| \le \tan\vphi\, |u_x|$, hence
$$
|\nabla u|^2 = u_x^2 + u_y^2 \le (1 + \tan^2 \vphi) u_x^2 = \frac{1}{\cos^2 \vphi}\, u_x^2,
$$
and
$$
\int_{\Om_2} u_x^2\, dx dy \ge \cos^2 \vphi \int_{\Om_2} |\nabla u|^2\, dx dy = \cos^2 \vphi\, I_2
\quad \Longrightarrow \quad \int_{\Om_2} u_y^2\, dx dy \le \sin^2 \vphi\, I_2.
$$

Fix $\vphi$ (and therefore, $r$) and let $\ve \to 0$; the integral $I_2$ goes to infinity, and $I_1$ is bounded, and therefore, $I_{1x} = \int_{\Om_1} u_x^2\, dx dy$ and $I_{1y} = \int_{\Om_1} u_y^2\, dx dy$ are also bounded. It follows that
$$
\frac{\int_\Om u_x^2\, dx\, dy}{\int_\Om u_y^2\, dx\, dy} =
\frac{\int_{\Om_1} u_x^2\, dx\, dy + \int_{\Om_2} u_x^2\, dx\, dy}{\int_{\Om_1} u_y^2\, dx\, dy + \int_{\Om_2} u_y^2\, dx\, dy} \ge
\frac{I_{1x} + \cos^2 \vphi\, I_2}{I_{1y} + \sin^2 \vphi\, I_2},
$$
hence the lower partial limit of the ratio satisfies
$$
\liminf_{\ve\to\infty}\frac{\int_\Om u_x^2\, dx\, dy}{\int_\Om u_y^2\, dx\, dy} \ge \cot^2 \vphi.
$$
Since $\vphi$ can be made arbitrarily small, the limit of the ratio is $+\infty$. This proves claim (a) of Theorem \ref{t1}.

\section{Proof of Theorem \ref{t1}: claim (b)}\label{st1b}

Here we give two proofs of claim (b).

{\it \underline{First proof}}
\vspace{2mm}

The intersections of vertical lines of support with $\Om$ are points. They will be called the {\it left corner} and the {\it right corner}.

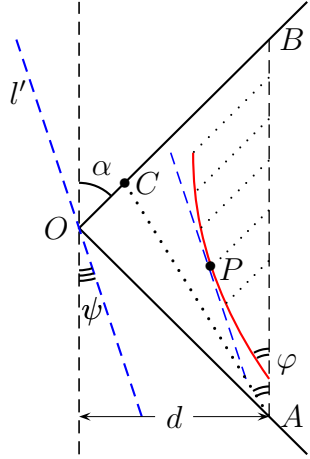
\begin{figure}[h]
\begin{picture}(0,170)
\scalebox{1}{

\rput(5,2.9){
\psline[linestyle=dashed,linewidth=0.6pt,linecolor=blue](2.195,-2)(1.205,1)
\psline[linestyle=dashed,linewidth=0.6pt](0,-3)(0,3)  \psline(3,-3)(0,0)(3,3)  \psline[linestyle=dashed,linewidth=0.5pt](2.5,-2.5)(2.5,2.5)
\rput(0.3,0.8){$\al$}  \psarc(0,0){0.6}{45}{90}
\psecurve[linecolor=red](3.35,-3)(2.5,-2)(1.93,-1)(1.6,0)(1.5,1)(1.6,2) \psarc(2.5,-2){0.4}{90}{120} \psarc(2.5,-2){0.32}{90}{120} \rput(2.73,-1.8){$\vphi$}
\psarc(2.5,-2.5){0.32}{90}{122} \psarc(2.5,-2.5){0.4}{90}{122}
\psline[linestyle=dotted](1.5,1)(2.5,2) \psline[linestyle=dotted](1.55,0.5)(2.5,1.45) \psline[linestyle=dotted](1.6,0)(2.5,0.9)
\psline[linestyle=dotted](1.76,-0.5)(2.5,0.24) \psline[linestyle=dotted](1.93,-1)(2.5,-0.43)
\rput(2.8,-2.5){$A$} \rput(2.8,2.5){$B$} \rput(-0.3,0){$O$}
\psline[linestyle=dashed,linewidth=0.8pt,linecolor=blue](0.825,-2.5)(-0.825,2.5) \rput(-0.8,1.8){$l'$}
\psarc(0,0){0.63}{270}{288} \psarc(0,0){0.7}{270}{288} \rput(0.17,-1.1){$\psi$} \psarc(0,0){0.56}{270}{288} \rput(0.17,-1.1){$\psi$}
\psdots(1.73,-0.5)(0.6,0.6) \rput(2,-0.55){$P$}
\psline[linecolor=black,linewidth=0.3pt,arrows=<-,arrowscale=1.5](0,-2.5)(1,-2.5)
\psline[linecolor=black,linewidth=0.3pt,arrows=->,arrowscale=1.5](1.5,-2.5)(2.5,-2.5) \rput(1.25,-2.5){$d$}
\psline[linestyle=dotted,linewidth=1.2pt](2.5,-2.5)(0.6,0.6) \rput(0.9,0.6){$C$}
}}
\end{picture}
\caption{The domain $\Om$ is contained in the angle $AOB$. The part of the level curve $u(x,y) = t$ contained in in $\Om_{++}$ passes through $P$.}
\label{fig:curve}
\end{figure}

For a certain $0 < \al \le \pi/4$, the domain $\Om$ is contained between two rays with the vertices at the left corner and the director vectors $(\sin\al, \cos\al)$ and $(\sin\al, -\cos\al)$ \, (the rays $OA$ and $OB$ in Fig.~\ref{fig:curve}).

Consider a nonzero function $u \in \UUU_{\Om}$. Without loss of generality assume that $\max u = 1$. Fix a value $0 < \vphi < \al$ and denote
$$
\Om_{u,\vphi} = \{ (x,y) \in \Om :\, |u_y| \le \tan\vphi\, |u_x| \}.
$$
Observe that
\beq\label{uxuy}
\int_\Om u_y^2\, dx\, dy \ge \int_{\Om\setminus\Om_{u,\vphi}} u_y^2\, dx\, dy \ge \tan^2 \vphi \int_{\Om\setminus\Om_{u,\vphi}} u_x^2\, dx\, dy.
\eeq

We are going to prove that there exist values $c_1$ and $c_2$ depending only on $\Om$ and $\vphi$ such that
\beq\label{a}
 \int_{\Om_{u,\vphi}} |\nabla u|^2\, dx\, dy \le c_1
\eeq
 and
\beq\label{b}
 \int_\Om u_y^2\, dx\, dy \ge c_2.
\eeq
It will then follow from \eqref{a}, \eqref{b}, and \eqref{uxuy} that
$$
\frac{\int_\Om u_x^2\, dx\, dy}{\int_\Om u_y^2\, dx\, dy}
\le \frac{\int_{\Om_{u,\vphi}} |\nabla u|^2\, dx\, dy}{\int_\Om u_y^2\, dx\, dy}
+ \frac{\int_{\Om\setminus\Om_{u,\vphi}} u_x^2\, dx\, dy}{\int_\Om u_y^2\, dx\, dy}
\le \frac{c_1}{c_2} + \cot^2 \vphi,
$$
and the proof of claim (b) of Theorem \ref{t1} will be finished.
 \vspace{2mm}

\underline{Proof of \eqref{a}}. \ We have $\Om_{u,\vphi} = \Om_{--} \cup \Om_{+-} \cup \Om_{-+} \cup \Om_{++}$, where
 $$\Om_{\ve_1\ve_2} =
 \Om_{\ve_1\ve_2}^{u,\vphi} = \{ (x,y) \in \Om :\, \ve_1 u_x \ge 0, \ \ve_2 u_y \ge 0, \ |u_y| \le \tan\vphi\, |u_x| \}, \quad \text{with} \ \ve_1,\, \ve_2 \in \{ +, -\}.
 $$
It suffices to prove that there exists $c > 0$ depending only on $\Om$ and $\vphi$ such that
$$
\int_{\Om_{++}} |\nabla u|^2\, dx\, dy \le c;
$$
the proofs for the domains $\Om_{+-}$,\, $\Om_{-+}$,\, $\Om_{--}$ are similar.

   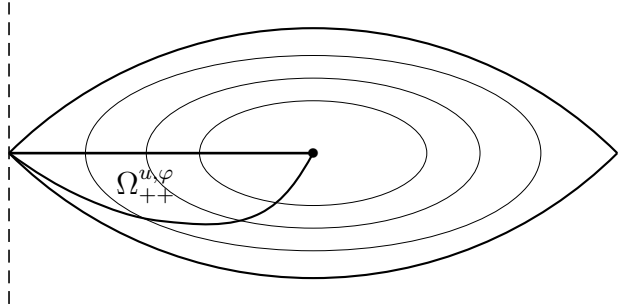
\begin{figure}[h]
\begin{picture}(0,125)
\scalebox{1}{

\rput(7.5,1.8){
\psarc(0,-4){5.657}{45}{135} \psarc(0,4){5.657}{225}{315}
\psdots(0,0)
\psellipse[linewidth=0.3pt](0,0)(3,1.3)  \psellipse[linewidth=0.3pt](0,0)(2.2,1) \psellipse[linewidth=0.3pt](0,0)(1.5,0.7)
\psline[linewidth=1pt](-4,0)(0,0)
\rput(-2.2,-0.45){\scalebox{1}{$\Om_{++}^{u,\vphi}$}}
\psecurve(-5,0.9)(-4,0)(-3,-0.6)(-2,-0.9)(-0.7,-0.8)(0,0)(0.3,1)  
\psline[linestyle=dashed,linewidth=0.6pt](-4,-2)(-4,2)

}}
\end{picture}
\caption{Several level curves of $u$ and a set $\Om_{++}^{u,\vphi}$ are shown.}
\label{fig:level}
\end{figure}

Consider a system of coordinates $t,\, s$ in $\Om$, where $t = u(x,y)$ and $s$ is a natural parameter along level curves of $u$; see Fig.~\ref{fig:level}. One has
\beq\label{dif}
dt\, ds = |\nabla u|\, dx\, dy.
\eeq

We are going to obtain estimates from above for $|\nabla u|$ and for the lengths of level curves $u(x,y) = t$ in $\Om_{++}$.

Take a level curve $u(x,y) = t$ in $\Om_{++}$ and denote by $d = d(t)$ the maximum distance between a point of the curve and the left vertical line of support. (Note that the maximum is attained at the lower point of the curve.) Let $O$ be the left corner and $AB$ be the vertical segment at the distance $d$ from the left vertical support line with the endpoints on the two rays with the vertex at $O$ and the director vectors $(\sin\al, \cos\al)$ and $(\sin\al, -\cos\al)$; see Fig.~\ref{fig:curve}. The length $l = l(t)$ of the level curve satisfies the inequality
\beq\label{l}
l \le |AB| = 2d \cot\al.
\eeq

Indeed, draw the straight lines parallel to the upper ray through the points of the segment $AB$ (dotted lines in Fig.~\ref{fig:curve}). They form the angle $\al$ with the segment and the angles $\psi + \al$ with the level curve, with $0 \le \psi \le \vphi$. The condition $0 < \vphi < \al \le \pi/2 - \al$ implies that the angles formed with the level curve lie in the interval $[\al, \pi/2].$ This means that the length of the arc of the curve between two infinitesimally close dotted lines is not greater than the length of the part of the segment $AB$ contained between the same lines. Summing up over all infinitesimal arcs, one concludes that $l \le |AB|$.

Let us now estimate the modulus of the gradient at a regular point, say $P$, of this curve. Denote by $\psi$,\, $0 \le \psi \le \vphi$, the angle between the tangent line to the curve at $P$ and the vertical. Draw the line $l'$ through $O$ parallel to the tangent line and denote by $h(P)$ the distance between these lines. Now consider the tangent plane to the graph of $u$ at the point $(P, t)$. The intersection of this plane with the horizontal plane $\RRR^2 \times \{0\}$ is a line parallel to the tangent line (and therefore, to $l'$) and disjoint with the interior of $\Om$. It follows that the distance between the line of intersection and the tangent line through $P$ is $\ge h(P)$. This implies that
$$|\nabla u(P)| \le t/h(P).$$

Draw the line through $A$ with the director vector $(-\sin\vphi, \cos\vphi)$ (dotted line in Fig.~\ref{fig:curve}), and let $C$ be the point of its intersection with $OB$. Since $P$ is contained in the triangle $ACB$, the distance $h(P)$ between $P$ and $l'$ is greater than or equal to the distance between $\bigtriangleup ACB$ and $l'$. Using geometric argument, one easily finds that the infimum of the latter distance over $\psi$ is attained at $\psi = 0$ and is equal to $d\, \dfrac{\sin(\al-\vphi)}{\sin(\al+\vphi)}$. Thus,
$$
h(P) \ge d\, \dfrac{\sin(\al-\vphi)}{\sin(\al+\vphi)},
$$
and therefore, the gradient of $u$ at $P$ satisfies
\beq\label{nabla}
|\nabla u| \le \dfrac{t\,\sin(\al+\vphi)}{d\,\sin(\al-\vphi)}.
\eeq
Using \eqref{dif}, \eqref{l}, and \eqref{nabla}, one obtains
$$
\int_{\Om_{++}} |\nabla u|^2\, dx\, dy = \int_0^1 \int_0^{l(t)}  |\nabla u|\, ds\, dt
$$ $$
\le \int_0^1 \Big( \int_0^{l(t)} \dfrac{t\,\sin(\al+\vphi)}{d(t)\,\sin(\al-\vphi)}  \, ds \Big) dt
\le \int_0^1 \dfrac{2t\, \sin(\al+\vphi)}{\sin(\al-\vphi)}\, \cot\al\, dt = \dfrac{\sin(\al+\vphi)}{\sin(\al-\vphi)}\, \cot\al.
$$
Inequality \eqref{a} is proved.
\vspace{2mm}

\underline{Proof of \eqref{b}}. \ Let the function $u$ attain its maximum at $\xi_0 = (x_0,y_0)$, that is, $u(x_0,y_0) = 1$. Let $\tilde u$ be the smallest function from $\UUU_\Om$ satisfying $\tilde u(x_0,y_0) = 1$. Clearly, the graph of $\tilde u$ is composed of line segments with one endpoint at $(\xi_0,1)$ and the other one on $\pl\Om \times \{0\}$, and $\tilde u \le u$.

Theorem 1 jointly with generalization 3 of Section 2 in the paper \cite{LRP} imply that if a function $f : \RRR^2 \to \RRR$ is convex and two convex functions $u_1$ and $u_2$ satisfy $u_2 \le u_1$ in  $\Om$ and $u_1 = u_2$ on $\pl\Om$, then $\int_\Om f(\nabla u_2)\, dx\, dy \ge \int_\Om f(\nabla u_1)\, dx\, dy.$ Now taking the convex functions $u_1 = -\tilde u$ and $u_2 = -u$ and using that the function $g(\z) = f(-\z)$ is also convex, one gets
$$
\int_\Om g(\nabla u)\, dx\, dy \ge \int_\Om g(\nabla\tilde u)\, dx\, dy.
$$
Taking $g(\z_1,\z_2) = \z_2^2$, one obtains
\beq\label{ineq11}
\int_\Om u_y^2\, dx\, dy \ge \int_\Om \tilde u_y^2\, dx\, dy.
\eeq

Let $h$ be the maximum distance between a point of $\Om$ and a line of support to $\Om$. The modulus of gradient at each regular point of $\tilde u$ satisfies $|\nabla\tilde u| \ge 1/h.$

Recall that a point $\xi \in \pl\Om$ is regular if and only if the line of support at $\xi$ is unique, and in this case the line is called the {\it tangent line} at $\xi$. Since both vertical lines of support are angular, there exists $0 < \bt < \pi/2$ such that each tangent line forms an angle $\ge \bt$ with the vertical.

Let a point $\xi \in \pl\Om$ be regular, and let the tangent line at $\xi$ form an angle $\vphi$ with the vertical, $\bt \le \vphi \le \pi/2$. Then each point $(x,y)$ on the open interval with the endpoints $\xi_0$ and $\xi$ is a regular point of $\tilde u$, and $|\tilde u_y| = \sin\vphi\, |\nabla\tilde u| \ge \sin\bt\, |\nabla\tilde u|$. It follows that
\beq\label{ineq22}
\int_\Om \tilde u_y^2\, dx\, dy \ge \frac{\sin^2 \bt}{h^2}\, |\Om|.
\eeq
Inequalities \eqref{ineq11} and \eqref{ineq22} imply \eqref{b} with $c_2 = \dfrac{\sin^2 \bt}{h^2}\, |\Om|$.
\vspace{2mm}

{\it \underline{Second proof}}
\vspace{2mm}

We need to prove that there is a constant $M$ such that, for every continuous concave
function~$u$ with bounded gradient,
which vanishes on the boundary of~$\Omega$, we have
 $\, \int_{\Omega} u_x^2 \, dx\, dy \, \le \, M \,
\int_{\Omega} u_y^2 \, dx\, dy \,$.

Without loss of generality it may be assumed
that~$\Omega$ has the extreme left point~$O$ at the origin  and  the projection of $\Om$
 to the $x$-axis is the segment~$[0,2]$, i.e., its extreme right point
 has coordinates~$(2, c)$.

Furthermore, after the linear change of variables
$x \, = \, \tilde x\, , \ y =  \frac{c}{2}\, \tilde x \, + \, \tilde y$ the domain~$(x, y) \in \Omega$
is mapped to a domain~$(\tilde x, \tilde y) \in \tilde \Omega $  with the extreme right point~$D = (2,0)$. The Jacobian of this transform is equal to one
and~$u_{y}\, = \, u_{\tilde y}, \ u_{x}\, = \,
 u_{\tilde x} \, -  \, \frac{c}{2}\, u_{\tilde y}$. Therefore,
$$
\int_{\Omega} u_y^2  \, dx\, dy \ = \
\int_{\tilde \Omega} u_{\tilde y}^2  \, d\tilde x\, d\tilde y \ ;
\qquad
\int_{ \Omega} u_x^2  \, dx\, dy
 \ = \
\int_{\tilde \Omega} \Bigl(  u_{\tilde x} \,  - \, \frac{c}{2}\,u_{\tilde y}\Bigr)^2 \,
d\tilde x\, d\tilde y \, .
$$
The inequality~$(a+b)^2 \, \le \, 2a^2 + 2b^2$ implies that
$$
\int_{ \Omega} u_x^2  \, dx\, dy
 \ \le \
2\, \int_{\tilde \Omega} \, u_{\tilde x}^2 \,
d\tilde x\, d\tilde y \ + \
2\, \int_{\tilde \Omega} \Bigl(\frac{c}{2}\, u_{\tilde y}\Bigr)^2  \,
d\tilde x\, d\tilde y \ .
$$
Therefore, if the ratio~$\dfrac{\int_{\tilde \Omega}
u_{\tilde x}^2\, d\tilde x d\tilde y}{\int_{\tilde \Omega}
u_{\tilde y}^2\, d\tilde x d\tilde y}$ does not exceed a constant~$\tilde M$, then
$\dfrac{\int_{\Omega}
u_{ x}^2\, d x \,d y}{\int_{ \Omega}
u_{ y}^2\, d x \, d x}$ does not exceed $\, M \, = \, 2\tilde M \, + \, \frac{c^2}{2}$. Thus,
it suffices to prove the theorem for the domain~$\tilde \Omega$. To simplify the notation
we will  assume that the original domain~$\Omega$ has the extreme right point~$D = (2,0)$,
so we shall not use tildes.

Let~$\Omega$ be bounded below and above by the graphs of functions~$\varphi_1(x)$
and~$\varphi_2(x)$ respectively. So, $\varphi_1$ is convex,~$\varphi_2$
is concave on the segment~$[0,2]$ and both those functions
 vanish at its ends.
\smallskip

It suffices to prove the inequality for the
the left part~$\frac12 \, \Omega \, = \, \bigl\{ (x, y) \in \Omega: \, 0\le y\le 1 \bigr\}$, then the same argument can be applied to the right part. Denote
$$
I_x \ = \ \int_{0}^1\, \, dx \, \int_{\varphi_1(x)}^{\varphi_2(x)} u_x^2\, dy \ ; \qquad
I_y \ = \ \int_{0}^1\, \, dx \, \int_{\varphi_1(x)}^{\varphi_2(x)}  u_y^2\, dy
$$
We need to show that~$I_x \, \le \, M\, I_y$.
\smallskip

\noindent \textbf{Lemma~1}. {\em At every point
$(x,y) \, \in \, \frac12 \, \Omega$, we have $\, \bigl|u_x\bigl| \, \le \,
\frac1x\, \bigl(\, u\, - \, y\, u_y \bigr)$.}
\smallskip

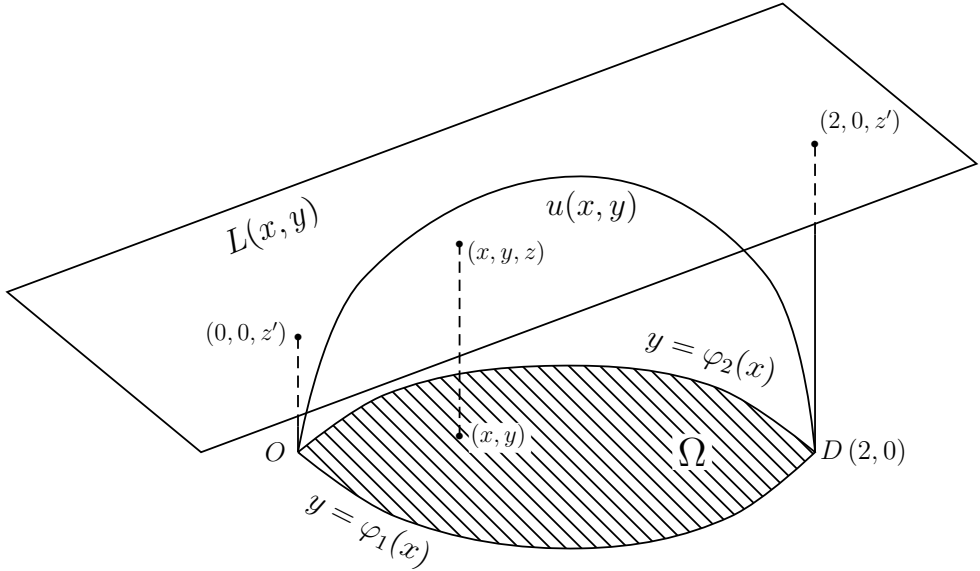
\begin{figure}[h]
\begin{picture}(0,200)
\scalebox{0.85}{

\rput(5.5,1.5){

\pscustom[fillstyle=vlines,fillcolor=yellow]{
\pscurve(0,0)(1.5,-1)(4.4,-1.5)(6.7,-1)(8,0)
\pscurve[linestyle=dashed](8,0)(6.5,1)(4,1.35)(1.35,0.9)(0,0)
}

\pspolygon[linewidth=0pt,linecolor=white,fillstyle=solid,fillcolor=white](2.65,0)(3.6,0)(3.6,0.45)(2.65,0.45)

\pspolygon[linewidth=0pt,linecolor=white,fillstyle=solid,fillcolor=white](5.8,-0.27)(5.8,0.28)(6.4,0.28)(6.4,-0.27)

\psdots[dotsize=3pt](2.5,0.25)(2.5,3.25) (0,1.8)(8,4.8)
\pspolygon(-1.5,0)(-4.5,2.5)(7.5,7)(10.5,4.5) \pscurve(0,0)(1,2.7)(4.5,4.3)(7.25,2.75)(8,0)
\psline[linestyle=dashed](2.5,0.25)(2.5,3.25)
     \psline(0,0)(0,0.4)       \psline[linestyle=dashed](0,0.4)(0,1.8)
     \rput(-0.8,1.85){\scalebox{0.9}{$(0,0,z')$}}     \rput(8.7,5.15){\scalebox{0.9}{$(2,0,z')$}}
         \psline(8,0)(8,3.4)       \psline[linestyle=dashed](8,3.4)(8,4.8)
\rput(6.1,0){\scalebox{1.5}{$\Om$}} \rput(3.1,0.25){\scalebox{0.9}{$(x,y)$}} \rput(8.75,0){\scalebox{1}{$D\, (2,0)$}}  \rput(-0.35,0){\scalebox{1}{$O$}} \rput(3.22,3.1){\scalebox{0.9}{$(x,y,z)$}}  \rput{23}(-0.4,3.5){\scalebox{1.25}{$L(x,y)$}}
\rput(4.55,3.82){\scalebox{1.2}{$u(x,y)$}}
\rput{-27}(1.1,-1.2){\scalebox{1.17}{$y = \vphi_1(x)$}}    \rput{-16}(6.4,1.5){\scalebox{1.17}{$y = \vphi_2(x)$}}
}}
\end{picture}
\caption{The tangent plane to the graph of $u$.}
\label{fig:plane}
\end{figure}

\noindent\textbf{Proof}. Denote by~$L(x,y)$ the tangent plane to the
graph of the function~$u$ at a regular point $(x,y)\in \Omega$; see Fig.~\ref{fig:plane}. The equation of this
plane in~$\re^3$ is
$$
(x', y', z') \in L(x,y) \qquad \Leftrightarrow \qquad \ z'-z \ = \ (x'-x)\, u_x(x,y) \ + \ (y'-y)\, u_y(x,y)\, ,
$$
where $z\, =\, u(x,y)$.
 Since $u$ is concave, it follows that
 $L(x,y)$ is located  above
the graph. i.e., $z' \, \ge \, u(x',y')$.
At the point~$O =  (0,0)$, we have~$u = 0$ and therefore,
$z'\ge  0$.  Substituting the point~$(0,0,z')$ to the equation of the plane, we get
$z + (-x )u_x + (-y)u_y \, = \, z' \ge 0$. Thus,
$\, u  \, -\, x \, u_x \, - \, y\, u_y \, \ge \, 0$ and so,
$u_x \, \le \,
\frac1x\, \bigl(\, u\, - \, y\, u_y \bigr)$.\
 This proves the lemma in the case~$u_x(x,y) \ge 0$.

If $u_x(x,y) < 0$, then we
apply the same argument to the
point~$D  = (2,0)$, where  we also have~$z' \ge 0$.
Substituting the point~$(2,0,z')$ to the equation of~$L(x,y)$, we obtain
 $z + (2-x )u_x + (-y)u_y \, \ge 0$,
which implies $u_x \, \ge \,
\frac{1}{2-x}\, \bigl(\, - u\, + \, y\, u_y \bigr)$.
Consequently, $|u_x| \, = \, -u_x \, \le \, \frac{1}{2-x}\,
\bigl(\, u\, - \, y\, u_y \bigr)$. The right hand side must be positive,
and  since~$x\le 1$, we have
$2-x \ge x$, which implies  $|u_x| \,  \le \, \frac{1}{x}\,
\bigl(\, u\, - \, y\, u_y \bigr)$. This completes the proof of the lemma.

{\hfill $\Box$}
\smallskip

 Now we continue the proof of the inequality. Applying Lemma~1 and the inequality $(a+b)^2 \le 2a^2 + 2b^2$, we get~$\  I_x \ = $
 $$
 \int_{0}^{1}\, \, dx\,  \int_{\varphi_1(x)}^{\varphi_2(x)} \, u_x^2 \, dy \ \le \
 \int_{0}^{1}\, \frac{dx}{x^2}\, \int_{\varphi_1(x)}^{\varphi_2(x)} \,
 \bigl(\, u\, - \, y\, u_y \bigr)^2 \, dy \ \le \  \int_{0}^{1}\, \frac{dx}{x^2}\, \int_{\varphi_1(x)}^{\varphi_2(x)} \,
\bigl( 2u^2\, + \, 2y^2\, u_y^2 \bigr)\, dy  \ =
$$
\begin{equation}\label{eq.ineq}
  2\, \int_{0}^{1}\, \frac{dx}{x^2}\, \int_{\varphi_1(x)}^{\varphi_2(x)} \,
u^2\, \, dy \ + \
  2\,  \int_{0}^{1}\, dx\,  \int_{\varphi_1(x)}^{\varphi_2(x)} \,
\left(\frac{y}{x} \right)^2\, u_y^2 \, dy \, .
\end{equation}
The convex domain~$\Omega$ lies
between two its tangents at the point~$O$.
Consequently, for every~$(x,y) \in \Omega$, we have $\varphi_1'(0) \le \dfrac{y}{x} \le \varphi_2'(0)$.
In the notation~$m = \max \{|\varphi_1'(0)|\, , \, |\varphi_2'(0)|\}$, we get
$\left(\dfrac{y}{x} \right)^2 \, \le m^2$ and therefore, the
 second integral in the sum~(\ref{eq.ineq}) does not
 exceed\linebreak $\int_{0}^{1}\, dx\, \int_{\varphi_1(x)}^{\varphi_1(x)} \, m^2 \,
 \, u_y^2 \, dy \, = \, m^2\, I_y$. To estimate the first
 integral in~(\ref{eq.ineq}), we use the Wirninger inequality (see, e.g., Section 7.7 of the book \cite{HR}): for an absolutely continuous function~$v$
 on~$[0,1]$ with a derivative from~$L_2$ and with~$v(0) = v(1) = 0$, we have
$$
\int_{0}^1  v^2(t) \, dt \ \le \
\frac{1}{\pi^2}\int_{0}^1  \dot v^2(t) \, dt
$$
(the equality is attained for $v(t) = \sin(\pi t)$).
After the change of variables~$y= (1-t)\varphi_1(x) \, + \, t\varphi_2(x)$, we obtain
$$
\int_{\varphi_1(x)}^{\varphi_2(x)}  v^2(y) \, dy \quad \le \quad
\left(  \frac{\varphi_2(x) - \varphi_1(x)}{\pi} \right)^2
\int_{\varphi_1(x)}^{\varphi_2(x)}  \dot v^2(y) \, dy\, .
$$
Now substitute this inequality to the first integral in~(\ref{eq.ineq}):
$$
 \int_{0}^{1}\, \frac{dx}{x^2}\, \int_{\varphi_1(x)}^{\varphi_2(x)} \,
u^2\, \, dy \quad \le \quad  \int_{0}^{1}\,
\, \frac{dx}{\pi^2}\, \left(  \frac{\varphi_2(x) - \varphi_1(x)}{x} \right)^2
  \, \int_{\varphi_1(x)}^{\varphi_2(x)} \,
u_y^2\, \, dy
$$
Note that~$\dfrac{\varphi_2(x) - \varphi_1(x)}{ x} \, \le \,
\varphi_2'(0) - \varphi_1'(0)\, \le \, 2m$.
Therefore, the latter integral does not exceed
$$
\frac{4m^2}{\pi^2}\,  \int_{0}^{1}\, dx\,\int_{\varphi_1(x)}^{\varphi_2(x)} \,
u_y^2\, dy \quad = \quad \frac{4m^2}{\pi^2} \, I_y\, .
$$
Combining those estimates, we conclude that
$\, I_x \, \le \, 2\, m^2 \, \bigl(1 + \frac{4}{\pi^2}
\bigr)\, I_y$,
which completes the proof.

\bigskip

\section{Proof of Theorem \ref{t2}}\label{st2}

Take an ortonormal coordinate system $x,\, y$ so as the eigendirections of the form $f''(0,0)$ coincide with the $x$-axis and the $y$-axis. In this coordinate system the Taylor decomposition of $f$ at the origin up to the quadratic term takes the form $f(\z) = f(0,0) + \al\z_1 + \bt\z_2 + \lam_1\z_1^2 + \lam_2\z_2^2 + o(|\z|^2) \quad \text{as} \ \, |\z| \to 0, \quad \text{where} \quad \z = (\z_1, \z_2).$ The contribution of the linear part $f(0,0) + \al\z_1 + \bt\z_2$ to the integral is a constant, $\int_\Om (f(0,0) + \al u_x + \bt u_y)\, dx\, dy = f(0,0) |\Om|$, therefore one can assume without loss of generality that $f(0,0) = \al = \bt = 0$, and write
\beq\label{f}
f(\z) = \lam_1\z_1^2 + \lam_2\z_2^2 + o(|\z|^2) \qquad \text{as} \ \, |\z| \to 0.
\eeq
Recall that the $C^1$ norm of $u$ is defined by \eqref{C1norm}. Note that at any regular point of $\Om$ one has
\beq\label{nab}
|\nabla u| = \sqrt{u_x^2 + u_y^2} \le |u_x| + |u_y| \le \|u\|_{C^1}.
\eeq

\vspace{2mm}

\underline{Proof of claim (a).} \ Let $f''(0,0)$ be positive definite; then $\lam_1 > 0$ and $\lam_2 > 0$. Choose $0 < \ve < \min \{ \lam_1,\, \lam_2 \}$. There exists $\del > 0$ such that if $|\z| < \del$ then the function $o(|\z|^2)$ in \eqref{f} satisfies $|o(|\z|^2)| \le \ve|\z|^2$. It follows that if $\|u\|_{C^1} < \del$ then, by \eqref{nab}, at any regular point of $\Om$ one has $|\nabla u| < \del$, hence $f(\nabla u) \ge (\lam_1-\ve)u_x^2 + (\lam_2-\ve)u_y^2 \ge 0$.
This implies that $F(u) \ge 0$ for all $u$ in the $\del$-neighborhood of zero, that is, $u \equiv 0$ is a local minimum of $F$.
\vspace{2mm}

\underline{Proof of claim (b).} \ Let $f''(0,0)$ be negative definite, and so, $\lam_1 < 0$ and $\lam_2 < 0$. Similarly to (a), choose $0 < \ve < \min \{ -\lam_1,\, -\lam_2 \}$. Again, there exists $\del > 0$ such that if $|\z| < \del$ then $o(|\z|^2)$ in \eqref{f} satisfies $|o(|\z|^2)| \le \ve|\z|^2$. Thus, if $\|u\|_{C^1} < \del$ then at any regular point of $\Om$ one has $|\nabla u| < \del$, hence $f(\nabla u) < (\lam_1+\ve)u_x^2 + (\lam_2+\ve)u_y^2 \le 0$. If, additionally, $u \not\equiv 0$ and hence, $|\nabla u| \ne 0$ at a certain point, then $f(\nabla u) < 0$ at this point. This implies that $F(u) < 0$ for all $u \not\equiv 0$ in the $\del$-neighborhood of the zero function, that is, $u \equiv 0$ is not a local minimum of $F$.
\vspace{2mm}

\underline{Proof of claim (c).} \ Let $f''(0,0)$ be indefinite; then $\lam_1$ and $\lam_2$ have different signs. Let, say, $\lam_1 = -a$ and $\lam_2 = b$ with $a$ and $b$ positive. Denote
$$
F_\ve(u) = \int_\Om \big( (-a+\ve)u_x^2 + (b+\ve)u_y^2 \big)\, dx\, dy.
$$
The lines of support parallel to the positive eigendirection are vertical.

(i) Suppose that at least one vertical line of support is not angular. Take $0 < \ve < a$. By claim (a) of Theorem \ref{t1}, there is a nonzero $\tilde u \in \UUU_\Om$ such that
 $$
 \int_\Om \tilde u_{x}^2\, dx\, dy > \dfrac{b+\ve}{a-\ve} \int_\Om \tilde u_{y}^2\, dx\, dy,
 $$
and so, $F_\ve(\tilde u) < 0$. There exists $\del > 0$ such that if $|\z| < \del$ then the function $o(|\z|^2)$ in \eqref{f} satisfies $|o(|\z|^2)| \le \ve|\z|^2$. Thus, if $\|u\|_{C^1} < \del$ then $f(\nabla u) \le (-a+\ve)u_x^2 + (b+\ve)u_y^2$ at any regular point of $\Om$, and therefore, $F(u) \le F_\ve(u)$. It follows that for $0 < t < \dfrac{\del}{\|\tilde u\|_{C^1}}$,\, $F(t \tilde u) \le F_\ve(t \tilde u) = t^2 F_\ve(\tilde u) < 0$. Thus, $u \equiv 0$ is not a local minimum of $F$.

(ii) Suppose now that both vertical lines of support are angular. Then by claim (b) of Theorem \ref{t1},
\beq\label{supK}
\sup_{\stackrel{u\in\UUU_\Om}{u\not\equiv0}}\, \frac{\int_\Om u_x^2\, dx\, dy}{\int_\Om u_y^2\, dx\, dy} = K,
\eeq
where $K = K(\Om) > 0$ is a real value.

Let $b/a < K$. Choose $\ve > 0$ sufficiently small such that
\beq\label{veK}
\frac{b+\ve}{a-\ve} < K.
\eeq
Repeating the above argument, one concludes that for an appropriate $\del > 0$, if $\|u\|_{C^1} < \del$ then $f(\nabla u) \le (-a+\ve)u_x^2 + (b+\ve)u_y^2$ at any regular point, and hence, $F(u) \le F_\ve(u)$. In view of \eqref{supK} and \eqref{veK}, there exists $\tilde u \in \UUU_\Om$ such that
 $$
 \int_\Om \tilde u_{x}^2\, dx\, dy > \dfrac{b+\ve}{a-\ve} \int_\Om \tilde u_{y}^2\, dx\, dy,
 $$
and so, $F_\ve(\tilde u) < 0$. It follows that for $0 < t < \dfrac{\del}{\|\tilde u\|_{C^1}}$,\, $F(t \tilde u) \le F_\ve(t \tilde u) = t^2 F(\tilde u) < 0$. Thus, $u \equiv 0$ is not a local minimum of $F$.

Let now $b/a > K$. Choose $\ve > 0$ sufficiently small such that $(b-\ve)/(a+\ve) > K$. By \eqref{supK}, all nonzero functions $u \in \UUU_\Om$ satisfy
$$
 \int_\Om u_{x}^2\, dx\, dy < \dfrac{b-\ve}{a+\ve} \int_\Om u_{y}^2\, dx\, dy.
$$
whence $F_{-\ve}(u) > 0$. For an appropriate $\del > 0$, if $\|u\|_{C^1} < \del$ and $u \not\equiv 0$ then $f(\nabla u) \ge (-a-\ve)u_x^2 + (b-\ve)u_y^2$, hence $F(u) \ge F_{-\ve}(u) > 0$. It follows that $u \equiv 0$ is a local minimum of $F$.

\section*{Acknowledgements}

The work of A. Plakhov was supported by CIDMA through FCT within projects UIDB/04106/2020 (https://doi.org/10.54499/UIDB/04106/2020), UIDP/04106/2020 (https://doi.org/10.54499/UIDP/04106/2020), and 2022.03091.PTDC\\ (https://doi.org/10.54499/2022.03091.PTDC). The research of V.Protasov was performed with the support of the Theoretical Physics and Mathematics Advancement Foundation “BASIS”.

\end{document}